\documentclass[twoside,12pt]{article}
\usepackage{amsmath}
\usepackage{enumitem}
\usepackage{amssymb}  
\newtheorem{theorem}{Theorem}[section]
\newtheorem{lemma}{Lemma}[section]
\newtheorem{proposition}{Proposition}[section]
\newtheorem{definition}{Definition}[section]

\def\Sum{\displaystyle\sum}

\def\Min{\displaystyle\min}

\def\Sup{\displaystyle\sup}
\def\Max{\displaystyle\max}
\def\Inf{\displaystyle\inf}
\def\Lim{\displaystyle\lim}

\DeclareMathOperator*{\esssup}{ess\,sup}
\begin{document}

\begin{center}
{\bf \large On a Generalized Prandtl-Batchelor Free Boundary Problem with a Singularity on a Stratified Lie Group   }\vskip 0.5cm {\bf By} \vskip
1cm {\bf Sabri BENSID }\\
\
\\ {\em  Address} : Laboratoire d'Analyse Nonlinéaire et Mathématiques Appliquées,\\ Department of Mathematics,\\ Faculty of Sciences,\\ University of Tlemcen,\\ B.P.119, Tlemcen 13000,
Algeria.\\  Mail: $ edp\_sabri@yahoo.fr$

\end{center}%
\begin{abstract}

We investigate a class of elliptic free boundary problems, including a  generalized Prandtl–Batchelor type problem with a singularity on a stratified Lie group. The associated energy functional is nondifferentiable, which precludes the direct application of standard variational techniques. To address this challenge, we approximate the energy functional using a sequence of $C^1$-functionals and establish the existence of at least two distinct nontrivial solutions. Our analysis combines elliptic regularity theory with the mountain pass theorem in a non-Euclidean setting.
\end{abstract}

\section{Introduction}
 During
the last few decades, the study of elliptic free boundary problems has been thoroughly investigated within the Euclidean framework. This type of problem arises in various physical fields, such as fluid dynamics. One interesting problem studied in this context is the Prandtl–Batchelor free boundary problem, initially studied in  $2$ dimensions by Prandtl and Batchelor. See \cite{Batchelor1, Batchelor2} for more details.\\\\
Recently, in \cite{perrera class},  K. Perera studied a more general form, as follows:
\begin{equation}\label{1}
  \left\{\begin{array}{ll}  -\Delta u=\lambda \chi_{\{u>1\}}g(x,(u-1)_+) & \quad \mbox{
in }\ \Omega\setminus F(u),\\[0.3cm]| \nabla u^+|^2-| \nabla u^-|^2=2 &\quad \mbox{ in }\
  F(u),
\\[0.3cm]u =0 &\quad \mbox{ on }\
  \partial \Omega,
\end{array}
\right.
 \end{equation}
 where $\Omega $ is a bounded domain in $\mathbb{R}^N,$ $ N\geq 2$ with a smooth boundary, $\lambda>0,$ $(u-1)_+=\Max(u-1,0),$ $F(u)=\partial \{u>1\}$ is the free boundary of $u$ and $\chi_{\{u>1\}}$ is the characteristic
function of the set. Also, $\nabla u^{\pm}$ are the limit of $\nabla u$ from the sets $\{u>1\}$ and $\{u\leq 1\}^{\circ}$ respectively and $g$ is a locally H\"{o}lder continuous function satisfying some assumptions.\\\\
Using the important ingredient  of the uniform Lipschitz continuity result of Caffarelli et
al \cite{caffa-annal} and an interesting approximation, the author proves the existence of two distinct nontrivial solutions. See \cite{perrera class} and the references therein. Some related results concerning problem $(1)$ are given in \cite{Caflish,Elcrat}.\\
The special case $g(x,s)=1$ is the well-known Prandlt-Batchelor free
boundary problem. When $g(x,s)=0,$ we refer the reader to \cite{alt-caffarelli,alt-caffarelli-fried,caffa-annal,weiss1,weiss2}. Other  generalizations like $g(x,s)=s^{p-1},$ $p>2$ are given in  \cite{Frie Plasma, Fri Liu, Temam, Temam1}.\\\\
More recently, considerable attention has been given to PDEs defined on stratified Lie groups (see the definition in Section 2). The challenge lies in determining whether results that are well-established for Euclidean domains (in a commutative setting) remain valid in the context of a noncommutative framework, particularly for free boundary problems. In this sense, in \cite{bensid-Lie},  we have studied  problem $(1)$ in the non-Euclidean setup. \\\\ More precisely, we consider the following problem
\begin{equation}\label{2}
  \left\{\begin{array}{ll}  -\mathcal{L}u=\lambda \chi_{\{u>1\}}g(x,(u-1)_+) & \quad \mbox{
in }\ \Omega\setminus F(u),\\[0.3cm]| \nabla_\mathbb{G}u^+|^2-| \nabla_\mathbb{G}u^-|^2=2 &\quad \mbox{ in }\
  F(u),
\\[0.3cm]u =0 &\quad \mbox{ on }\
  \partial \Omega.
\end{array}
\right.
 \end{equation}
The domain $\Omega\subset \mathbb{G}$  is  bounded where  $\mathbb{G}$ is a stratified Lie group and $\mathcal{L}$ is the sub-Laplacian which will be defined in the following Section. Finally, $\nabla_\mathbb{G}u^{\pm}$ are the limit of $\nabla_\mathbb{G}u$ from the sets $\{u>1\}$ and $\{u\leq 1\}^{\circ}$ respectively,
and $g:\Omega\times [0,\infty)\rightarrow [0,\infty)$ is a locally H\"{o}lder continuous function satisfying\\\\
$(g_1)$ Fore some $\alpha,\beta>0$ and $1<p<2,$
$$g(x,s)\leq \alpha+\beta s^{p-1},\quad \quad  \forall (x,s)\in \Omega\times[0,\infty).$$
$(g_2)\quad g(x,s)>0,\quad \forall x\in \Omega$ and $s>0.$\\\\
We have shown that there exists a positive solution $u$ satisfying the equation $ -\mathcal{L}u=\lambda \chi_{\{u>1\}}g(x,(u-1)_+)$ in the classical sense in $\Omega\setminus F(u)$ where the free boundary condition holds in the sense of viscosity and $u$ vanishes continuously on $\partial \Omega.$ See \cite{bensid-Lie} for further information.\\\\
One key result in this direction is the following monotonicity lemma, provided by Choudhuri and Repov$\breve{s}$ in \cite{Repsov}.
\begin{lemma}\cite{Repsov}
Let $u>0$ be a Lipschitz continuous function on the unit ball $B_1(0)\subset \mathbb{G}$ satisfying the distributional inequalities
$$\pm \mathcal{L}u\leq \left(\frac{\lambda}{\varepsilon}\chi_{\{|u-1|<\varepsilon\}}(x)\mathfrak{F}(|\nabla_\mathbb{G}u|)+A\right),$$
for constants $A>0,$ $0<\varepsilon<1.$ Suppose further that $\mathfrak{F}$ is a continuous function such that $\mathfrak{F}(t)=o(t^2)$ near infinity. Then, there exist $C=C(N,A)>0$ and $\int_{B_1(0)}u^2 dx$ but not on $\varepsilon$ such that
$$\esssup_{x\in B_{1/2}(0)}\{|\nabla_\mathbb{G}u|\}\leq C.$$
\end{lemma}
Using this result, Choudhuri and Repov$\breve{s}$  \cite{Repsov} prove the existence of  solutions to the following problem
\begin{equation}\label{2}
  \left\{\begin{array}{ll}  -\mathcal{L}u=\lambda (u-1)_+^2 f & \quad \mbox{
in }\ \Omega\setminus F(u),\\[0.3cm]| \nabla_\mathbb{G}u^+|^2-| \nabla_\mathbb{G}u^-|^2=2 &\quad \mbox{ in }\
  F(u),
\\[0.3cm]u =0 &\quad \mbox{ on }\
  \partial \Omega,
\end{array}
\right.
 \end{equation}
 where $f\in L^{\infty}(\Omega)$ is a positive bounded function.\\\\
In this article, we study the existence of solutions to the following problem
 \begin{equation}\label{3}
  \left\{\begin{array}{ll}  -\mathcal{L}u=\lambda \chi_{\{u>1\}}g(x,(u-1)_+)+\frac{\beta}{u^{\delta}} & \quad \mbox{
in }\ \Omega\setminus F(u),\\[0.3cm]| \nabla_\mathbb{G}u^+|^2-| \nabla_\mathbb{G}u^-|^2=2 &\quad \mbox{ in }\
  F(u),
\\[0.3cm]u =0 &\quad \mbox{ on }\
  \partial \Omega,
\end{array}
\right.
 \end{equation}
where $0<\delta<1$ and $\beta$ is a positive parameter.\\\\
Our study was motivated by the question of whether the results obtained in \cite{bensid-Lie} could be extended when a singular term is added. The distinguishing feature of this work, setting it apart from our paper \cite{bensid-Lie}  ,( Perera's approach \cite{perrera class} in the Euclidean setting), is the effective management of the singular term, which prevents the associated energy functional from being $C^1$ at $u=0$. This challenge is what makes it impossible to directly apply the results from the variational framework.\\\\
We define the variational functional associated to problem $(\ref{3})$  by
$$E(u)=\int_{\Omega}\left[\frac{|\nabla_{\mathbb{G}}u|^2}{2}+\chi_{\{u>1\}}(x)-\lambda G(x,(u-1)_+)\right]dx-\frac{\beta}{1-\delta}\int_{\Omega}(u^+)^{1-\delta}dx,$$
where $G(x,s)=\int_{0}^{s}g(x,t)dt,\quad s\geq 0.$\\
So, the functional \( E \) is not of class \( C^1 \) due to the term \( \int (u^+)^{1-\gamma} \, dx \). Furthermore, it is nondifferentiable because of the term \( \int \chi_{\{u > 1\}}(x) \, dx \).\\\\
Note also that the analogue of problem $(\ref{3})$ in the Euclidean  setting was studied by Choudhuri and Repov$\breve{s}$ in \cite{Repsov-sing} where the authors used  the existence of a solution to the following classical problem
 \begin{equation}\label{5}
  \left\{\begin{array}{ll}  -\Delta u=\frac{f(x)}{u^{\delta}} & \quad \mbox{
in }\ \Omega,
\\[0.3cm]u =0 &\quad \mbox{ on }\
  \partial \Omega,
\end{array}
\right.
 \end{equation}
 where $f>0,$  of class $C^{\alpha},$ $0<\alpha<1$ and $\delta>0.$\\
 Problem $(\ref{5})$ was considered by Lazer and McKenna  \cite{Lazer}, where they proved the existence and uniqueness of a solution $u\in C^{2,\alpha}(\Omega)\cap C(\overline{\Omega}).$\\\\
Recently, using the fixed point theorem argument,  The authors in \cite{chen} proved the existence of solutions to the following problem involving the sub-Laplacian on the Heisenberg group
\begin{equation}\label{6}
  \left\{\begin{array}{ll}  -\Delta_{H^N} u= \frac{f(x)}{u^{\delta}} & \quad \mbox{
in }\ \Omega,
\\[0.3cm]u =0 &\quad \mbox{ on }\
  \partial \Omega,
\end{array}
\right.
 \end{equation}
where $\Omega$ is a bounded domain in the Heisenberg group $H^N$ and $\Delta_{H^N}$ is  usually called the Kohn Laplacian or Heisenberg Laplacian. See also \cite{wang} for subelliptic problem  on the Heisenberg group
with mixed singular and power type nonlinearity.\\
The analogue of problem $(\ref{6})$ with the p-sub-Laplacian ($p>1$) was investigated in \cite{Garain} in the
setting of stratified Lie groups for $0<\delta<1.$ For an extension of problem $(\ref{6}),$ the reader can also see the work of  Sahu et al \cite{Sahu}.\\\\
Very recently, there has been growing interest in nonlocal subelliptic operators, particularly the fractional sub-Laplacian on stratified Lie group. Interested readers may refer to \cite{Ghosh0,Ghosh,Ghosh2} and the references therein.\\\\
Now, we define the Sobolev space $W^{1,p}(\Omega)$ where $\Omega$ is an open subset on a stratified Lie group as
$$W^{1,p}(\Omega)=\{u\in L^p(\Omega):\quad |\nabla_\mathbb{G}u|\in L^p(\Omega)\}.$$
A norm on this space
$$\|u\|_{1,p}=\|u\|_p+\|\nabla_\mathbb{G} u\|_{p}.$$
However, the Sobolev space $W_0^{1,p}(\Omega)$ is defined as
$$W_0^{1,p}(\Omega)=\{u\in W^{1,p}(\Omega),\quad u=0\quad \hbox{on}\quad \partial \Omega\},$$
where $u=0$ on $\partial \Omega$ is the usual trace sense. ( See section $2$ for more details).\\
By a solution, we mean a function $u\in W_0^{1,2}(\Omega), u>0$ a.e in $\Omega$ satisfying
$$\int_{\Omega}\nabla_{\mathbb{G}}u\nabla_{\mathbb{G}}\varphi dx=\lambda \int_{\Omega}g(x,(u-1)_+)\varphi dx+\beta\int_{\Omega}u^{-\delta}\varphi dx,\quad \forall \varphi \in W_0^{1,2}(\Omega).$$
We denote by $\lambda_1$ the first eigenvalue of $(-\mathcal{L})$ and the corresponding eigenfunction by $\varphi_1.$ By testing problem $(\ref{3})$ with $\varphi_1,$ the following weak formulation must be satisfied for $u$ to be a weak solution of problem
\begin{equation}\label{31}
 \lambda_1\int_{\Omega}u\varphi_1 dx=\int_{\Omega}\nabla_{\mathbb{G}}u\nabla_{\mathbb{G}}\varphi_1 dx=\lambda \int_{\Omega}g(x,(u-1)_+)\varphi_1 dx+\beta\int_{\Omega}(u^+)^{-\delta}\varphi_1 dx.
\end{equation}
Thus, there exists a constant, denoted by $\beta^*$, depending on $\lambda$, such that the inequality $\beta t^{–\delta} + \lambda g(x, (t – 1)_+) > \lambda_1 t$ holds for all $t > 0.$ This contradicts $(\ref{31})$. Therefore, $0 < \beta < \beta^*$.\\\\
So, the first result of this paper is the following theorem

\begin{theorem}
 Let  $0 < \beta < \beta^*$. If we denote by $u_{\beta}$ the solution of the following problem
\begin{equation}\label{sing}
  \left\{\begin{array}{ll}  -\mathcal{L}u=\frac{\beta}{u^{\delta}} & \quad \mbox{
in }\ \Omega,\\[0.3cm] u>0 &\quad \mbox{ in }\
  \Omega,
\\[0.3cm]u =0 &\quad \mbox{ on }\
  \partial \Omega,
\end{array}
\right.
 \end{equation}
 then $u_{\beta}<u$ a.e in $\Omega,$ where $u$ is the weak solution of problem $(\ref{3}).$\\
 Moreover, for $u\in W_{loc}^{1,2}(\Omega),$ the Radon measure $\mathcal{M}=\mathcal{L}u+\beta u^{-\delta}$ is nonnegative and supported on $\Omega\cap\{u<1\} $ for $\beta\in (0,\beta^*).$
\end{theorem}
In the next result, we will give the existence of solutions of problem  $(\ref{3})$ with some properties.
\begin{theorem}
Suppose that $(g_1)$ and $(g_2)$ are satisfied. Then, there exists $\lambda^*,\beta^*>0$ such that for all $\lambda>\lambda^*,$\quad $0<\beta<\beta^*,$ problem $(\ref{3})$ has two Lipschitz continuous solutions say $u_0,u_1\in W_0^{1,2}(\Omega)\cap C^2(\overline{\Omega}\setminus F(u))$ that satisfy the equation $-\mathcal{L}u=\lambda \chi_{\{u>1\}}g(x,(u-1)_+)+\beta u^{-\delta}$ in the classical sense in $\Omega\setminus F(u),$ the free boundary condition in the
generalized sense, and vanish continuously on $\partial \Omega.$\\
Moreover, \\
1) $E(u_0)<-\mathcal{H}(\Omega)\leq -\mathcal{H}(\{u=1\})<E(u_1),$ where $\mathcal{H}(.)$ denotes the Haar measure, and hence $u_0,u_1$ are nontrivial and distinct.\\
2) $0<u_1\leq u_0,$ the regions $\{u_0<1\}\subset \{u_1<1\}$ are connected and the regions $\{u_0>1\}\supset \{u_1>1\}$ are nonempty.\\
3) $u_0$ is a minimizer of $E$ but $u_1$ is not a minimizer of $E.$
\end{theorem}
The rest of this paper is organized as follows.
\newpage
\section{Preliminaries: stratified Lie groups and Sobolev spaces}
In this section, we provide the necessary tools in order to study problem $(\ref{1}).$ The reader who is familiar with these notions may go directly to the next section. We begin by giving the definitions of homogeneous Lie group.
\begin{definition}
Let $\mathbb{G}=(\mathbb{R}^N,*)$ be a Lie group on $\mathbb{R}^N.$ We say that $G$ is a homogeneous if any $\delta>0,$ there exists an automorphism $T_{\delta}:\mathbb{G}\rightarrow \mathbb{G}$ defined by
$$T_{\delta}(x)=(\delta^{r_1}x_1,\delta^{r_2}x_2,...,\delta^{r_N}x_N),\quad r_i>0,i=1,...,N\quad\hbox{for}\quad x=(x_1,...,x_N). $$
The map $T_{\delta}$ is called a dilation on $\mathbb{G}.$ The positive number $Q=\Sum_{i=1}^{N}r_i$ called homogeneous dimension of the group $\mathbb{G}.$
\end{definition}
\begin{definition}
We say that $\mathbb{G}=(\mathbb{R}^N,*)$ is a homogeneous Carnot group or a homogeneous stratified Lie group if the following properties hold\\
1) $\mathbb{R}^N$ can be split as $\mathbb{R}^N=\mathbb{R}^{N_1}\times...\mathbb{R}^{N_r}$ for some natural number $N_1,N_2,...,N_r$ such that $N_1+N_2+...+N_r=N$ and the dilation $T_{\delta}:\mathbb{R}^N\rightarrow \mathbb{R}^N,$
$$T_{\delta}(x)=(\delta^{1}x^{(1)},\delta^{2}x^{(2)},...,\delta^{r}x^{(r)}),\quad x^{(i)\in \mathbb{R}^N}$$
is an automorphism of the group $\mathbb{G}$ for every $\delta>0.$\\
2) If $N_1$ is as above, let $Z_1,...,Z_{N_1}$ be the left invariant vector fields on $\mathbb{G}$ such that $Z_i(0)=\frac{\partial}{\partial x_i}\large|_0$ for $i=1,...,N,$ then
\begin{equation}\label{5}Rank(Lie\{Z_1,...,Z_{N_1}\})(x)=N,\quad \forall x\in \mathbb{R}^N
\end{equation}
The vector fields $Z_1,...,Z_{N_1}$ are called the (Jacobian) generators of $\mathbb{G}$ whereas any basis for $span{Z_1,...,Z_{N_1}}$ is called a systm of generator of $\mathbb{G}$. We also say that $\mathbb{G}$ has step $r$ and $N_1$ generators.\\
Note also that $(\ref{5})$ is the well-known the H\"{o}rmander condition.
\end{definition}
An interesting example is the $\mathbb{H}^1=(\mathbb{R}^3,*) $ known as the Heisenberg-Weyl group on $\mathbb{R}^3$ which is a stratified Lie group of step $2$ and $2$ generators. Indeed, it is a homogeneous Lie group with dilation
$$T_{\delta}(x_1,x_2,x_3)=(\delta x_1,\delta x_2,\delta^2 x_3)$$
and for $$Z_1=\frac{\partial}{\partial x_1}+2x_2 \frac{\partial }{\partial x_3},\quad Z_2=\frac{\partial}{\partial x_2}-2 x_1\frac{\partial}{\partial x_3},$$
we have
$$rank(Lie\{Z_1,Z_2\})(x)=3,\quad \hbox{for every}\quad x\in \mathbb{R}^3.$$
Thus, $1)$ and $2)$ of the above definition are satisfied.\\\\
In the sequel, we let $N_1=N$ and let $Z=(Z_1,...,Z_N)$ be any family of vector fields on $\mathbb{R}^N.$
\begin{definition}
A Lipschitz continuous curve $\gamma:[0,T]\rightarrow \mathbb{R}^N, T\geq 0$ is $Z-$admissible if there exists a bounded functions $h=(h_1,h_2,...,h_N):[0,T]\rightarrow\mathbb{R}^N$ such that
$$\gamma'(t)=\Sum_{j=1}^N h_j Z_j(\gamma(t)),\quad \hbox{for a.e}\quad t\in [0,T].$$
In general, the function $h$ is not unique.
\end{definition}
Then, the following definition makes sense.
\begin{definition}(Carnot-Carathéodory distance)
For any $x,y\in \mathbb{G},$ the Carnot-Carathéodory distance is defined by
$$d_Z(x,y)=inf\{T\geq 0, \hbox{then there exits an admissible}\quad \gamma:[0,T]\rightarrow \mathbb{G}\ \hbox{such that}\quad \gamma(0)=x,\gamma(T)=y\}.$$
If the above set is empty put $d_Z(x,y)=0.$\\
If $d_Z(x,y)<+\infty$ for all $x,y\in \mathbb{G},$ then $(\mathbb{G},d_Z))$ is a metric space known as the Carnot-Carathéodory space.
\end{definition}
Now, we define the Sub-Laplacian ( or horizontal Laplacian) on $\mathbb{G}$ as
$$\mathcal{L}:=Z_1^2+...+Z_N^2.$$
The horizontal gradient on $\mathbb{G}$ is defined as
$$\nabla_\mathbb{G}:=(Z_1,...,Z_N),$$
and the horizontal divergence on $\mathbb{G}$ is defined by
$$div_\mathbb{G}v=\nabla_\mathbb{G}.v.$$
Thus, the sub-Laplacian $\mathcal{L}$ can be written as
$$\mathcal{L}u=div_\mathbb{G}(\nabla_\mathbb{G}u).$$
On the other hand, for $x\in \mathbb{G},$ the Haar measure ( denoted by $dx$) is the Lebesgue measure on $\mathbb{R}^N.$ So, for $1\leq p\leq \infty,$ the space $L^p(\mathbb{G})$ or simply $L^p$ denote the usual Lebesgue space on $\mathbb{G}$ with respect to the Haar measure with the  norm
$$||.||_{L^p}=\|.\|_{L^p(\mathbb{G})}=\|.\|_p$$
given  for $p\in [1,+\infty)$ by
$$\|f\|_p=\left(\int_\mathbb{G}|f(x)|^p dx\right)^{1/p},$$
and for $p=\infty$ by
$$\|f\|_{\infty}=\Sup_{x\in \mathbb{G}}|f(x)|.$$
Here, the supremum with respect to the Haar measure. We recall that the Haar measure for a measurable subset of $\mathbb{G}$ is denoted by $\mu(.)$ and satisfies the properties
$$\mu(T_{\delta}(M))=\delta^Q\mu(M)$$
where $M\subset \mathbb{G}$ is a measurable set.\\\\
For the variational setting of our problem, recall that the Sobolev space $W^{1,p}(\Omega)$ where $\Omega$ is an open subset on a stratified Lie group is defined as
$$W^{1,p}(\Omega)=\{u\in L^p(\Omega):\quad |\nabla_\mathbb{G}u|\in L^p(\Omega)\}.$$
A norm on this space
$$\|u\|_{1,p}=\|u\|_p+\|\nabla_\mathbb{G} u\|_{p}.$$
However, we have
$$W_0^{1,p}(\Omega)=\{u\in W^{1,p}(\Omega),\quad u=0\quad \hbox{on}\quad \partial \Omega\},$$
where $u=0$ on $\partial \Omega$ is the usual trace sense. This space is separable and uniformly convex Banach space.\\
If $\Omega$ is bounded, then for $1\leq p\leq Q,$ $W_0^{1,p}(\Omega)$ is continuously embedded in $L^q$ for $1\leq q\leq p^*=\frac{Qp}{Q-p}$ and compact if $1\leq q< p^*.$\\
Also, the Sobolev space $W_{loc}^{1,p}(\Omega)$ is the space of all functions in $W^{1,p}(K)$ for any compact set $K\subset \Omega.$\\
We state also the analogue of divergence theorem in the Euclidean setup for sub-Laplacian on stratified Lie group. This result is very important for the proofs of our results.
\begin{proposition}\cite{Ruzhansky}
Let $f_n\in C^1(\Omega)\cap C(\overline{\Omega}),$ $n=1,2,...,N.$ Then, for each $n,$ we have
$$\int_{\Omega}Z_n f_n d\nu=\int_{\partial \Omega}f_n<Z_n,d\nu>.$$
Consequently,
$$\int_{\Omega}\Sum_{n=1}^{N}Z_n f_n d\nu=\int_{\partial \Omega} \Sum_{n=1}^{N} f_n<Z_n,d\nu>.$$
\end{proposition}
Finally, we use in this paper the following theorem
\begin{theorem}\label{Am-Rabin}
Let $E$ be a $C^1-$functional defined on a Banach space $X$. Assume that $E$ satisfies the (PS) condition and there exist an open set $U\subset X,$ $u_0\in U$ and $u_1\in X\setminus \overline{U}$ such that
$$\Inf_{u\in \partial U}E(u)>\Max\{E(u_0),E(u_1)\}.$$
Then, $E$ has a critical point at the level
$$c:=\Inf_{\gamma\in \Gamma} \Max_{u\in \gamma([0,1])}E(u)\geq \Inf_{u\in \partial U}E(u),$$
where $\Gamma=\{\gamma\in C([0,1],X),\quad \gamma(0)=u_0,\gamma(1)=u_1\}$ is the class of paths in $X$ joining $u_0$ and $u_1.$
\end{theorem}
\section{Proof of Theorem 1.1}
Following \cite{Garain}, we assume the existence of $u_{\beta}$ solution of problem $\ref{sing}.$ So, for any $0<\beta<\beta^*,$ we have
$$0=\int_{\Omega}\nabla_{\mathbb{G}}u_{\beta}\nabla_{\mathbb{G}}\varphi dx-\beta \int_{\Omega}u_{\beta}^{-\delta}\varphi dx,\quad \forall \varphi \in W_0^{1,2}(\Omega). $$
Let $U=\{x\in \Omega,\quad u(x)<u_{\beta}(x)\}.$ Using the weak formulation of $u$ and testing by $\varphi=(u_{\beta}-u)_+,$ we obtain
$$0\leq \int_{\Omega}\nabla_{\mathbb{G}}(u_{\beta}-u)\nabla_{\mathbb{G}}(u_{\beta}-u)_+ dx$$
$$=-\lambda \int_{\Omega}\chi_{\{u>1\}}g(x,(u-1)_+)(u_{\beta}-u)_+dx+\beta\int_{\Omega}(u_{\beta}^{-\delta}-u^{-\delta})(u_{\beta}-u)_+dx\leq 0.$$
This implies that $\|u_{\beta}-u\|_{1,2}=0$ and hence $|U|=0.$\\
Note also that since $u_{\beta}$ and $u$ are continuous ( see the following section for existence and regularity of $u$ of problem $(\ref{3}).$ It follows that $U=\emptyset.$ So, $u_{\beta}\geq u$ in $\Omega.$\\\\
Now, let $\widetilde{U}=\{x\in \Omega,\quad u(x)=u_{\beta}(x)\}.$\\
Remark that $\widetilde{U}$ is a measurable set. Then, for any $\eta>0,$ there exists a closed subset $V$ of $\widetilde{U}$ such that $|\widetilde{U}\setminus V|\leq \eta.$\\
Assume that $|\widetilde{U}|>0$ and we define the following test function $\varphi$ such that
$$\varphi(x)=\left\{\begin{array}{ll} 1  & \quad \mbox{if }\ x\in V\\[0.3cm]
				0<\varphi<1        & \quad \mbox{if }\ x\in \widetilde{U}\setminus V\\[0.3cm]
				 0        & \quad \mbox{if }\ x\in \Omega\setminus \widetilde{U}
			\end{array}
			\right.	
		$$
The function $u$ is a weak solution of problem $(\ref{3})$ and satisfies
$$0=\int_{\Omega}-\mathcal{L}u\varphi dx-\beta \int_{V}u^{-\delta}dx-\beta\int_{\widetilde{U}\setminus V}u^{-\delta}\varphi dx$$
$$-\int_V g(x,(u-1)_+)dx-\int_{\widetilde{U}\setminus V}g(x,(u-1)_+)\varphi dx$$
$$=-\int_V g(x,(u-1)_+)dx-\int_{\widetilde{U}\setminus V}g(x,(u-1)_+)\varphi dx<0.$$
This is a contradiction. Thus, $|\widetilde{U}|=0$ and $\widetilde{U}=\emptyset.$ So, $u>u_{\beta}$ in $\Omega .$\\\\
Now, to prove the second part of theorem $1.1,$ we choose $\varepsilon>0,$ $\beta\in (0,\beta^*)$ and a test function $\varphi^2\chi_{\{u<1-\varepsilon\}}$ where $\varphi\in C_0^{\infty}(\Omega).$ We obtain
$$0=-\int_{\Omega}\nabla_{\mathbb{G}}u\nabla_{\mathbb{G}}(\varphi^2 \Min(u-1+\varepsilon,0))dx+\beta\int_{\Omega} u^{-\delta}\varphi^2\Min(u-1+\varepsilon,0)dx$$
$$=\int_{\Omega\cap \{u<1-\varepsilon\}}\nabla_{\mathbb{G}}u\nabla_{\mathbb{G}}(\varphi^2 (u-1+\varepsilon)dx+\beta\int_{\Omega\cap \{u<1-\varepsilon\} } u^{-\delta}\varphi^2(u-1+\varepsilon)dx$$
$$=\int_{\Omega\cap \{u<1-\varepsilon\}}|\nabla_{\mathbb{G}}u|^2\varphi^2 dx+2\int_{\Omega\cap \{u<1-\varepsilon\}}\varphi(u-1+\varepsilon)\nabla_{\mathbb{G}}u\nabla_{\mathbb{G}}(\varphi(u-1+\varepsilon))dx$$$$+\beta\int_{\Omega\cap \{u<1-\varepsilon\} } u^{-\delta}\varphi^2(u-1+\varepsilon)dx$$
The integration by parts gives
$$\int_{\Omega\cap \{u<1-\varepsilon\}}|\nabla_{\mathbb{G}}u|^2\varphi^2 dx=-2\int_{\Omega\cap \{u<1-\varepsilon\}}\varphi(u-1+\varepsilon)\nabla_{\mathbb{G}}u\nabla_{\mathbb{G}}(\varphi(u-1+\varepsilon))dx$$$$-\beta\int_{\Omega\cap \{u<1-\varepsilon\} } u^{-\delta}\varphi^2(u-1+\varepsilon)dx$$
$$\leq C \int_{\Omega}u^2|\nabla_{\mathbb{G}}\varphi|^2dx-\beta\int_{\Omega}u^{1-\delta}\varphi^2dx$$
$$\leq C \int_{\Omega}u^2|\nabla_{\mathbb{G}}\varphi|^2dx,\quad \hbox{for some}\quad C>1.$$
Since, $\int_{\Omega}|u|^2dx<\infty,$ by passing to the limit $\varepsilon\rightarrow 0,$ we conclude that $u\in W_{loc}^{1,2}(\Omega).$\\\\
Finally, for nonnegative $\xi\in C_0^{\infty}(\Omega),$ we have
$$-\int_{\Omega}\widetilde{\nabla}_{\mathbb{G}}\xi\widetilde{\nabla}_{\mathbb{G}}u+\beta\int_{\Omega}u^{-\delta}\xi dx$$
$$=\int_{\Omega\cap\{0<u<1-2\varepsilon\}}\left[\widetilde{\nabla}_{\mathbb{G}}(\xi \Max\{\Min(2-\frac{1-u}{\varepsilon},0)\})\widetilde{\nabla}_{\mathbb{G}}u+\beta u^{-\delta}\xi \right]dx$$
$$+\int_{\Omega\cap\{1-2\varepsilon<u<1-\varepsilon\}}\left[\widetilde{\nabla}_{\mathbb{G}}(\xi \Max\{\Min(2-\frac{1-u}{\varepsilon},0)\})\widetilde{\nabla}_{\mathbb{G}}u+\beta u^{-\delta}\xi \right]dx$$
$$+\int_{\Omega\cap\{1-\varepsilon<u<1\}}\left[\widetilde{\nabla}_{\mathbb{G}}(\xi \Max\{\Min(2-\frac{1-u}{\varepsilon},0)\})\widetilde{\nabla}_{\mathbb{G}}u+\beta u^{-\delta}\xi \right]dx$$
$$+\int_{\Omega\cap\{1<u\}}\left[\widetilde{\nabla}_{\mathbb{G}}(\xi \Max\{\Min(2-\frac{1-u}{\varepsilon},0)\})\widetilde{\nabla}_{\mathbb{G}}u+\beta u^{-\delta}\xi \right]dx$$
$$\geq \int_{\Omega\cap\{1-2\varepsilon<u<1-\varepsilon\}}\left[\left(2-\frac{1-u}{\varepsilon}\right)\widetilde{\nabla}_{\mathbb{G}}\xi\widetilde{\nabla}_{\mathbb{G}}u+\xi\frac{|\nabla_{\mathbb{G}}u|^2}{\varepsilon}+\beta u^{-\delta}\xi\right]dx\geq 0.$$
\newline
So, if $\varepsilon\rightarrow 0,$ then $\mathcal{L}(u-1)_{-}\geq 0$ in the distributional sense. Thus, there exists a Radon measure $\mathcal{M}$ such that $\mathcal{M}=\mathcal{L}u+\beta u^{-\delta}$  and supported on $\Omega\cap\{u<1\} $ for $\beta\in (0,\beta^*).$
\section{Approximate solution}
We will begin by addressing the singular term by introducing a cut-off function $\varphi_{\beta}$ defined as follows
$$ \varphi_{\beta}(u)=\left\{\begin{array}{ll} u^{-\delta}  & \quad \mbox{if }\ u>u_{\beta}\\[0.3cm]
				u_{\beta}^{-\delta}        & \quad \mbox{if }\ u\leq u_{\beta}
			\end{array}
			\right.	
		$$
where $u_{\beta}$ is the solution of problem  $(\ref{sing}).$\\\\
We approximate the functional energy $E$ associated to problem $(\ref{3})$  by $C^1-$functionals. So, let $\beta:\mathbb{R}\rightarrow [0,2]$ be a function such that\\
$(1)$ $\beta\in C^{\infty}(\mathbb{R}),$ support of $\beta$ lies in $[0,1]$ and it is positive in $(0,1).$\\
$(2)$ $\int_{0}^{1}\beta(s)ds=1.$\\\\
Then, let $B(s)=\int_{0}^{s}\beta(t)dt$ where $B:\mathbb{R}\rightarrow [0,1]$ is a smooth and nondecreasing function such that
$$B(s)=0, \quad \hbox{for}\quad s\leq 0,\quad 0<B(s)<1\quad \hbox{for}\quad 0<s<1\quad \hbox{and}\quad B(s)=1\quad \hbox{for}\quad s\geq 1.$$
Hence, for $\varepsilon>0,$ let
$$E_{\varepsilon}(u)=\int_{\Omega}\left[\frac{|\nabla_\mathbb{G}|^2}{2}+B(\frac{u-1}{\varepsilon})-\lambda G_{\varepsilon}(x,(u-1)_+)-\beta\Phi_{\beta}(u)\right]dx$$
where $G_{\varepsilon}(x,s)=\int_{0}^{s}g_{\varepsilon}(x,t)dt,\quad s\geq 0,$ \quad $g_{\varepsilon}(x,s)=B(\frac{s}{\varepsilon})g(x,s),$ $\Phi_{\beta}(u)=\int_{0}^{u}\varphi_{\beta}(t)dt.$\\\\
The function $E_{\varepsilon}$ is of class $C^1$ and its critical points are the solutions of the following problem
\begin{equation}\label{app}
  \left\{\begin{array}{ll}  -\mathcal{L}u=-\frac{1}{\varepsilon}\beta(\frac{u-1}{\varepsilon})+\lambda g_{\varepsilon}(x,(u-1)_+)+\beta\varphi_{\beta}(u) & \quad \mbox{
in }\ \Omega,
\\[0.3cm]u =0 &\quad \mbox{ on }\
  \partial \Omega.
\end{array}
\right.
 \end{equation}
Let $\varepsilon_j\searrow 0,$ let $u_j$ be a critical point of $E_{\varepsilon_j}.$ We prove the following convergence result
 \begin{lemma}
 Assume $g_1)$ and $g_2).$ If $(u_j)$ is bounded in $W_0^{1,2}(\Omega)\cap L^{\infty}(\Omega),$ then there exists a Lipschitz continuous function $u$ on $\overline{\Omega}$ such that
 \begin{enumerate}[label=(\roman*)]
        \item $u_j\rightarrow u$ uniformly over $\overline{\Omega}.$
        \item $u_j\rightarrow u$ locally in $C^1(\overline{\Omega}\setminus \{u=1\}).$
        \item $u_j\rightarrow u$ strongly in $W_0^{1,2}(\Omega).$
        \item $E(u)\leq \Lim \inf E_{\varepsilon_j}(u_j)\leq \Lim \sup E_{\varepsilon_j}(u_j)\leq E(u)+\mathcal{H}\{u=1\}.$
    \end{enumerate}
 Moreover, $u$ satisfies $-\mathcal{L}u=\lambda \chi_{\{u>1\}}g(x,(u-1)_+)+\beta u^{-\delta}$ in the classical sense in $\Omega\setminus F(u)$ and $u$ satisfies the free boundary condition in the generalized sense and vanishes continuously on $\Omega.$
 \end{lemma}
 {\bf Proof of Lemma 4.1} Assume that $0<\varepsilon_j<1.$ So, $u_j$ is a solution of
 \begin{equation}\label{appsol}
  \left\{\begin{array}{ll}  -\mathcal{L}u_j=-\frac{1}{\varepsilon}\beta(\frac{u_j-1}{\varepsilon})+\lambda g_{\varepsilon}(x,(u_j-1)_+)+\beta\varphi_{\beta}(u_j) & \quad \mbox{
in }\ \Omega,
\\[0.3cm]u_j =0 &\quad \mbox{ on }\
  \partial \Omega.
\end{array}
\right.
 \end{equation}
 Since $(u_j)$ is bounded in $L^{\infty}(\Omega),$ then $0\leq g_{\varepsilon_j}(x,(u_j-1)_+)\leq A_0$ for some constant $A_0>0.$\\\\
 Let $v_0$ be a solution of the following problem
 \begin{equation}\label{pb A}
  \left\{\begin{array}{ll}  -\mathcal{L} v_0=\lambda A_0+\beta u_{\beta}^{-\delta} & \quad \mbox{
in }\ \Omega,
\\[0.3cm]v_0 =0 &\quad \mbox{ on }\
  \partial \Omega.
\end{array}
\right.
 \end{equation}
 Since $\beta \geq 0,$ then $-\mathcal{L}u_j\leq \lambda A_0+\beta u_{\beta}^{-\delta}.$ \\ By the maximum principle, we have
 $$0\leq u_j(x)\leq v_0(x),\quad \quad \forall x\in \Omega.$$
 From theorem 1.1, we have that $u_j>u_{\beta}$ in $\Omega$ for all $\beta\in(0,\beta^*),$ and since,  $\{u_j\geq 1\}\subset \{v_0\geq 1\},$ it follows that $v_0$ gives a uniform lower bound, say $d_0$ on the distance from the set $\{u_j\geq 1\}$ to $\partial \Omega.$\\\\
 Also, since the sequence $(u_j)$ is bounded with respect to $C^{2,\alpha}-$norm, so, from standard boundary regularity theory, it has  a convergent subsequence in $C^2-$norm on $\frac{d_0}{2}$ neighborhood of $\partial \Omega.$\\\\
  Observe that $0\leq \beta\leq 2\chi_{(-1,1)}, $ then
 $$\pm \mathcal{L}u_j=\pm \frac{1}{\varepsilon_j}\beta(\frac{(u_j-1)_+}{\varepsilon_j})\pm \lambda g_{\varepsilon_j}(x,(u_j-1)_+)+\beta u_j^{-\delta}$$
 $$\leq \frac{2}{\varepsilon_j}\chi_{\{|u_j-1|<\varepsilon_j\}}(x)+\lambda A_0+\beta u_{\beta}^{-\delta}.$$
 Note that for any subset $M$ relatively compact of $\Omega,$ we have $u_{\beta}\geq C_M,$ for some $C_M>0.$ ( See \cite{Garain}). Hence, using lemma 1.1, then there exists $A>0$ such that
 $$\esssup_{x\in B_{r/2}(x_0)}\{|\nabla_\mathbb{G}u_j|\}\leq \frac{A}{r},\quad \hbox{for}\quad r>0,$$
 for which $B_r(x_0)\subset \Omega.$\\
 Because $(u_j)$ is a sequence of Lipschitz continuous functions, then we have
 $$\Sup_{x\in B_{r/2}(x_0)}\{|\nabla_\mathbb{G}u_j|\}\leq \frac{A}{r}.$$
 Therefore, the family $(u_j)$ is uniformly Lipschitz continuous on the compact subsets, say $K$ such that $d(K,\partial \Omega)\geq \frac{d_0}{2}.$\\
 The application of Ascoli-Arzela gives a subsequence namely also $(u_j)$ that converge uniformly to a Lipschitz continuous function $u$ in $\Omega$ and finally, the Banach-Alaoglu, we have $u_j\rightharpoonup u$ in $W_0^{1,2}.$\\
 Now, we show that $u$ satisfies
 \begin{equation}\label{et}
 -\mathcal{L}u=\lambda \chi_{\{u>1\}}g(x,(u-1)_+)+\beta u^{-\delta}\quad \hbox{in}\quad \{u\neq 1\}.
 \end{equation}
 Remark that three scenarios may arise concerning the position of $u_{\beta}$ with respect to $(0,1).$ More precisely, the following cases can occur
 $$0<u_{\beta}<1<u,\quad 0<u_{\beta}<u<1,\quad 1<u_{\beta}<u.$$
 Fortunately, none of the cases present significant mathematical challenges in establishing the convergence result.\\
 So, let $\xi\in C_0^{\infty}(\{u>1\}).$ Then, on the support of $\xi,$ $u\geq 1+2\eta$ for $\eta>0.$\\
For large $j,$ $\eta_j<\eta$ and  $|u_j-u|<\eta,$ we have $u_j\geq 1+\eta_j$ on the support of $\xi.$\\
So, testing $(\ref{appsol})$ with $\xi$ gives
$$\int_{\Omega}\widetilde{\nabla}u_j\widetilde{\nabla}\xi dx=\lambda\int_{\Omega}g_{\varepsilon_j}(x,u_j-1)\xi+\beta\int_{\Omega}u_j^{-\delta}\xi dx,$$
where $\widetilde{\nabla}\eta=\Sum_{k=1}^{N}X_k\eta X_k.$\\
Passing to the limit $j\rightarrow +\infty,$ we have
$$\int_{\Omega}\widetilde{\nabla}u\widetilde{\nabla}\xi dx=\lambda\int_{\Omega}g(x,u-1)\xi+\beta\int_{\Omega}u^{-\delta}\xi dx,$$
Hence $u$ is a distributional (and thus a classical) solution of
\begin{equation}\label{E1}
 -\mathcal{L}u=\lambda g(x,(u-1)_+)+\beta u^{-\delta}\quad \hbox{in}\quad \{u> 1\}.
\end{equation}
By the same argument, taking $\xi\in C_0^{\infty}(\{u<1\}),$ we find $\eta>0$ such that $u\leq 1-2\eta$ giving also $u_j<1-\eta.$\\
testing $(\ref{appsol})$ with any nonnegative $\xi$ and passing to the limit yields
$$-\mathcal{L}u=  \beta u^{-\delta}\quad \hbox{in}\quad \{u<1\}.$$
To know the nature of the set $\{u<1\}^{0},$ we test $(\ref{appsol})$ by any nonnegative function and we use the fact that $B<1$ and also passing to the limit to obtain
$$-\mathcal{L}u\leq \lambda g(x,(u-1)_+)+  \beta u^{-\delta}\quad \hbox{in the distributional sense.}$$
Following theorem 1.1, $\mathcal{M}=\mathcal{L}u+\beta u^{-\delta}$ is a positive Radon  measure supported on $\Omega \cap \{u<1\}.$\\
By the classical result of regularity ( see section 9.4 in \cite{gilbarg}), we obtain that $u\in W_{loc}^{2,2}(\{u\leq 1\}^o)$ and hence $\mathcal{M}$ is actually supported on $\Omega\cap \partial \{u<1\}\cap \{u>1\},$ So $u$ satisfies $\mathcal{L} u=\beta u^{-\delta}$ on the set $\{u\leq 1\}^o.$\\\\
To prove ii), we see that  $u_j\rightarrow u$ with respect to $u$ with respect to $C^2$ norm in a neighborhood of $\partial \Omega$ in $\overline{\Omega}.$ Then, for $\eta>0,$ we have $u\geq 1+2\eta$ in the set $U\subset\subset\{u>1\}.$\\
So, for large $j$ with $\eta_j<\eta,$ we have $|u_j-u|<\eta$ in $\Omega$ and hence $u_j\geq 1+\eta_j$ in $U.$\\
From $(\ref{app}),$ we have
$$-\mathcal{L} u_j=\lambda g(x,(u_j-1))\quad \hbox{in}\quad U.$$
Observe that $g$ is locally H\"{o}lder continuous and $u_j\rightarrow u$ uniformly, then $g(x,(u_j-1))\rightarrow g(x,(u-1))$ in $L^p(U),$ $p>1.$ Also, since $W^{2,2}(U)\hookrightarrow C^1(U),$ then, we have $u_j\rightarrow u$ in $C^1(U)$ which implies that $u_j\rightarrow u$ in $C^1(\{u>1\}).$\\
The same method can be applied to prove that $u_j\rightarrow u$ in $C^1(\{u<1\}).$\\\\
Now, we prove $iii).$ we remark that $u_j\rightharpoonup u$ in $W_0^{1,2}(\Omega)$ and by the weak lower semicontinuity of the norm $\|.\|$ where $\|u\|=\|\nabla_\mathbb{G} u\|_{p}.$ We have
$$\|u\|\leq \varliminf \|u_j\|.$$
Hence, to show iii), it suffices  to prove that $\varlimsup \|u_j\|\leq \|u\|.$\\
So, we multiply $(\ref{app})$ by $(u_j-1)$ and then integrate by parts. Using also the fact that $\tan(\frac{t}{\delta_j})\geq 0$ for all $t>0.$ Then, we obtain
\begin{align}\int_{\Omega}|\nabla_\mathbb{G}|^2dx & \leq \lambda \int_{\Omega}g(x,(u_j-1)_+)(u_j-1)_+dx-\int_{\partial \Omega}u_j<X_i,dn>ds+\beta \int_{\Omega}u_j^{-\delta}(u_j-1)_+ dx \nonumber\\
&\rightarrow \lambda \int_{\Omega}g(x,(u-1)_+)(u-1)_+dx-\int_{\partial \Omega}u<X_i,dn>ds+\beta \int_{\Omega}u^{-\delta}(u-1)_+ dx
\end{align}
where $n$ is the outward unit normal to $\partial \Omega.$\\
Now, fix $0<\varepsilon<1$ and testing the equation $(\ref{E1})$ with $\xi=(u-1-\varepsilon)_+$ yields
\begin{equation}\label{I1}
\int_{\{u>1+\varepsilon\}}|\nabla_\mathbb{G}u|^2dx=\int_{\Omega}\lambda g(x,(u-1)_+)(u-1-\varepsilon)_+ +\beta \int_{\Omega}u^{-\delta}(u-1-\varepsilon)_+dx.
\end{equation}
On the other hand, integrating $(u-1+\varepsilon)_{-}\Delta u=\beta \int_{\Omega}u^{-\delta}(u-1-\varepsilon)_-$ over $\Omega$ gives
\begin{equation}\label{I2}
\int_{\{u<1-\varepsilon\}}|\nabla_\mathbb{G}u|^2dx=-(1-\varepsilon)\int_{\partial\Omega}u<X_i,dn>ds+\beta \int_{\Omega}u^{-\delta}(u-1-\varepsilon)_- .
\end{equation}
Adding, $(\ref{I1})$ and $(\ref{I2})$ and letting $\varepsilon\rightarrow 0$ gives
$$\int_{\Omega}|\nabla_\mathbb{G}u|^2dx=\int_{\Omega}\lambda g(x,(u-1)_+)(u-1)_+dx-\int_{\partial \Omega}u<X_i,dn>ds+\beta \int_{\Omega}u^{-\delta}(u-1)_+.$$
Since $\int_{\{u=1\}}|\nabla_\mathbb{G}u|^2dx=0.$ So, using $(11),$ we have
$$\varlimsup \int_{\Omega}|\nabla_\mathbb{G}u_j|^2dx\leq \int_{\Omega}|\nabla_\mathbb{G}u|^2dx$$
as wanted.\\
For the last statement $iv),$ we have
$$E_{\varepsilon_j}(u_j)=\int_{\Omega}\left[\frac{1}{2}|\nabla_\mathbb{G}u_j|^2+B\left(\frac{u_j-1}{\varepsilon_j}\right)\chi_{\{u\neq 1\}}-\lambda G_{\varepsilon_j}(x,(u_j-1)_+)-\beta u_j^{-\delta}(u_j-1)_+\right]dx$$$$+\int_{\{u=1\}}B\left(\frac{u_j-1}{\varepsilon_j}\right)dx.$$
Since $B\left(\frac{u_j-1}{\varepsilon_j}\right)\chi_{\{u\neq 1\}}$ converges pointwise to $\chi_{\{u>1\}}$ and is bounded by $1$
 and $G_{\varepsilon_j}(x,(u_j-1)_+)$ converges to $G(x,(u-1)_+),$ then the first integral converges to $E(u).$\\
 On the other hand, observe that
 $$0\leq \int_{\{u=1\}}B\left(\frac{u_j-1}{\varepsilon_j}\right)dx\leq \mathcal{H}\{u=1\}. $$
 To conclude the proof of Lemma $4.1,$ we choose $\overrightarrow{\varphi}\in C_0^1(\Omega,G)$ such that $u\neq 1$ a.e on the support of $\overrightarrow{\varphi}.$ \\
 Multiplying by $\Sum_{k=1}^{N}\varphi_k X_k u_n$ the weak formulation of $(\ref{appsol})$ and integrate over the set $\{1-\varepsilon^-< u_n<1+\varepsilon^+\},$ we get
 \begin{align}\label{I3}
   \int_{\{1-\varepsilon^-< u_n<1+\varepsilon^+\}}\left[-\Delta_\mathbb{G}u_n+\frac{1}{\varepsilon_n}\beta(\frac{u_n-1}{\varepsilon_n})\right]\Sum_{k=1}^{N}\varphi_k X_k u_n dx & = \nonumber\\
     =\int_{\{1-\varepsilon^-< u_n<1+\varepsilon^+\}}(\lambda g_{\varepsilon_n}(x,(u_n-1)_+)+\beta u_n^{-\delta})\Sum_{k=1}^{N}\varphi_k X_k u_ndx.
 \end{align}
 The term of the left in $(\ref{I3})$ can be simplified as
\begin{align}
\nabla_\mathbb{G}\left(\frac{1}{2}|\nabla_\mathbb{G}u_n|^2\overrightarrow{\varphi}-\left(\Sum_{k=1}^{N}X_k u_n \varphi_k\right)|\nabla_\mathbb{G}u_n|\right)+\Sum_{k=1}^{N}\Sum_{l=1}^{N}X_l \varphi_k u_n X_k u_n \nonumber\\-\frac{1}{2}|\nabla_\mathbb{G}u_n|^2\nabla_\mathbb{G}\overrightarrow{\varphi}
+\Sum_{k=1}^{N}\varphi_k X_k B(\frac{u_n-1}{\varepsilon_n}).
\end{align}
We integer by part to obtain
\begin{align}\label{I4}
  \int_{\{1-\varepsilon^-< u_n<1+\varepsilon^+\}}\frac{1}{2}|\nabla_\mathbb{G}u_n|^2\Sum_{k=1}^{N}\varphi_k<X_k,dn>-(\Sum_{k=1}^{N}X_k u_n \varphi_k)\Sum_{l=1}^{N}X_l u_n<X_l,dn> \nonumber  \\
  + \int_{\{1-\varepsilon^-< u_n<1+\varepsilon^+\}}B(\frac{u_n-1}{\varepsilon_n})\Sum_{k=1}^{N}\varphi_k<X_k,dn> \nonumber  \\
 = \int_{\{1-\varepsilon^-< u_n<1+\varepsilon^+\}}(\frac{1}{2}|\nabla_\mathbb{G}|^2\Sum_{k=1}^{N} X_k\varphi_k-\Sum_{k=1}^{N}\Sum_{l=1}^{N}X_k\varphi_lX_l u_n X_k u_n)\nonumber \\
 + \int_{\{1-\varepsilon^-< u_n<1+\varepsilon^+\}}\left[B(\frac{u_n-1}{\varepsilon_n})\Sum_{k=1}^{N}X_k \varphi_k+(\lambda g_{\varepsilon_n}(x,(u_n-1)_+)+\beta u_n^{-\delta})\Sum_{k=1}^{N}X_k \varphi_k\right]dx.
\end{align}
The integral on the left converges to
\begin{align}
& \int_{\left\{u=1+\epsilon^{+}\right\} \cup\left\{u=1-\epsilon^{-}\right\}}\left[\frac{1}{2}\left|\nabla_{\mathbb{G}} u\right|^{2} \sum_{k=1}^{N} \varphi_{k}\left\langle X_{k}, d n\right\rangle-\left(\sum_{k=1}^{N} X_{k} u \varphi_{k}\right) \sum_{l=1}^{N} X_{l} u\left\langle X_{l}, d n\right\rangle\right. \\
+ & \left.\int_{\left\{u=1+\epsilon^{+}\right\}} \sum_{k=1}^{N} \varphi_{k}\left\langle X_{k}, d n\right\rangle\right] \nonumber \\
= & \int_{\left\{u=1+\epsilon^{+}\right\} \cup\left\{u=1-\epsilon^{-}\right\}}\left[\left(1-\frac{1}{2}\left|\nabla_{\mathbb{G}} u\right|^{2}\right) \sum_{k=1}^{N} \varphi_{k}\left\langle X_{k}, d n\right\rangle-\sum_{k \neq l ; 1 \leq k, l \leq N} \varphi_{k} X_{l} u X_{k} u\left\langle X_{l}, d n\right\rangle\right]  \\
= & \int_{\left\{u=1+\epsilon^{+}\right\}}\left[\left(1-\frac{1}{2}\left|\nabla_{\mathbb{G}} u\right|^{2}\right) \sum_{k=1}^{N} \varphi_{k}\left\langle X_{k}, d n\right\rangle-\sum_{k \neq l ; 1 \leq k, l \leq N} \varphi_{k} X_{l} u X_{k} u\left\langle X_{l}, d n\right\rangle\right]\nonumber
\end{align}
$$
\begin{aligned}
& -\int_{\left\{u=1-\epsilon^{-}\right\}}\left[\left(\frac{1}{2}\left|\nabla_{\mathbb{G}} u\right|^{2}\right) \sum_{k=1}^{N} \varphi_{k}\left\langle X_{k}, d n\right\rangle-\sum_{k \neq l ; 1 \leq k, l \leq N} \varphi_{k} X_{l} u X_{k} u\left\langle X_{l}, d n\right\rangle\right] \\
= & \int_{\left\{1-\epsilon^{-}<u<1+\epsilon^{+}\right\}}\left(\frac{1}{2}\left|\nabla_{\mathbb{G}} u\right|^{2} \sum_{k=1}^{N} X_{k} \varphi_{k}-\sum_{k=1}^{N} \sum_{l=1}^{N} X_{k} \varphi_{l} X_{l} u X_{k} u\right) d x \\
& \int_{\left\{1-\epsilon^{-}<u<1+\epsilon^{+}\right\}}\left[\sum_{k=1}^{N} X_{k} \varphi_{k}+(\lambda g(x,(u-1)_+)+\beta u^{-\delta})  \sum_{k=1}^{N} X_{k} \varphi_{k}\right] d x,
\end{aligned}
$$

as $n \rightarrow \infty$.\\
Observe that the normal vector at the point $P$ on the set $\left\{u=1+\epsilon^{+}\right\} \cup\left\{u=1-\epsilon^{-}\right\}$is $n= \pm \frac{\nabla_{\mathbb{G}} u(P)}{\left|\nabla_{\mathbb{G}} u(P)\right|}$. So, equation $(18)$ becomes when $\varepsilon\rightarrow 0$
$$0=\lim _{\epsilon \rightarrow 0} \int_{\left\{u=1+\epsilon^{+}\right\}}\left[\left(1-\frac{1}{2}\left|\nabla_{\mathbb{G}} u\right|^{2}\right) \sum_{k=1}^{N} \varphi_{k}\left\langle X_{k}, d n\right\rangle\right] \\
-\lim _{\epsilon \rightarrow 0} \int_{\left\{u=1-\epsilon^{-}\right\}}\left[\left(\frac{1}{2}\left|\nabla_{\mathbb{G}} u\right|^{2}\right) \sum_{k=1}^{N} \varphi_{k}\left\langle X_{k}, d n\right\rangle\right] .$$
Hence, we conclude that $u$ satisfies the free boundary condition in the sense of viscosity.
\section{Proof of Theorem 1.2}
In this section, we prove the existence of solutions of problem $(\ref{3}).$ First, we have the following lemma.
\begin{lemma}
Under the assumptions $g_1)$ and $g_2),$ there exists $\lambda^*,\beta^*>0$ such that for all $\lambda>\lambda^*,$ $0<\beta<\beta^*,$ $\varepsilon<\varepsilon_0(\lambda),$ the functional $E_{\varepsilon}$ has a minimizer $u_0^{\varepsilon}>0$ that satisfies
\begin{equation}\label{etoile}
  E_{\varepsilon}(u_0^{\varepsilon})\leq m_1(\lambda)+2\lambda \varepsilon a_0 \mathcal{H}(\Omega)<0,
\end{equation}
where $m_1(\lambda)=\Inf_{u\in W_0^{1,2}(\Omega)} E(u)$ and $\varepsilon_0(\lambda)=\Min\left(\frac{|m_1(\lambda)|}{2\lambda a_0 \mathcal{H}(\Omega)},\left(\frac{pa_0}{a_1}\right)^{\frac{1}{p-1}}\right).$
\end{lemma}
{\bf Proof of Lemma 5.1}\\ Under the condition $g_1),$ we have
$$E_{\varepsilon}(u)\geq \int_{\Omega}\left[\frac{1}{2}|\nabla_{\mathbb{G}}u|^2-\lambda\left(a_0(u-1)_+ +\frac{a_1 (u-1)_+}{p}\right)-\frac{\beta}{1-\delta}u^{1-\delta}\right]dx.$$
Since $1<p<2,$ then $E_{\varepsilon}$ is bounded from below and coercive. So, $E_{\varepsilon}$ satisfies the Palais-Smale condition and it possesses a minimizer $u_0^{\varepsilon}.$\\
Also, by the assumption $g_2),$ we have $G(x,t)\geq 0$ for all $x\in \Omega$ and $t>0.$ Thus, for any $u\in W_0^{1,2}(\Omega)$ with $u>1$ on a set of positive measure. Hence, $E(u)\rightarrow -\infty$ as $\lambda\rightarrow +\infty.$\\\\
It follows that there exists $\lambda^*>0$ such that $\forall \lambda >\lambda^*,$ we have \begin{equation}\label{neg}
                                                                                             m_1(\lambda)=\Inf_{u\in W_0^{1,2}(\Omega)} E(u)<-\mathcal{H}(\Omega).
                                                                                           \end{equation}
On the other hand, because $B(\frac{t-1}{\varepsilon})\leq \chi_{(1,\infty)}(t),$ for all $t,$ then
$$E_{\varepsilon}(u)-E(u)\leq \lambda\int_{\Omega}\left[G(x,(u-1)_+)-G_{\varepsilon}(x,(u-1)_+)\right]dx$$
$$=\lambda\int_{\Omega}\int_{0}^{(u-1)_+}\left[\left(1-B(\frac{t}{\varepsilon})\right)g(x,t)\right]dt dx$$
$$\leq \lambda\int_{\Omega}\int_{0}^{\varepsilon}g(x,t)dt dx$$
$$\leq \lambda(a_0\varepsilon+\frac{a_1}{p}\varepsilon^p)\mathcal{H}(\Omega).$$
For $\varepsilon<\varepsilon_0(\lambda),$ we obtain $(\ref{etoile}).$ So, $E_{\varepsilon}(u_0^{\varepsilon})\leq 0=E_{\varepsilon}(0).$
Then, $u_0^{\varepsilon}$ is a nontrivial solution of problem $(\ref{appsol}).$ $\Box$\\\\
Next, we will prove that $E_{\varepsilon}$ admits a second nontrivial critical point denoted by $u_1^{\varepsilon}.$
\begin{lemma}
For $\lambda>\lambda^*,$ $0<\beta<\beta^*,$  there exists $C(\lambda)$ such that $\forall \varepsilon<\varepsilon_0(\lambda),$ the functional $E_{\varepsilon}$ has a second critical point $0<u_1^{\varepsilon}<u_0^{\varepsilon}$ that verifies
$$C(\lambda)\leq E_{\varepsilon}(u_1^{\varepsilon})\leq \frac{1}{2}\|u_0^{\varepsilon}\|^2+\mathcal{H}(\Omega).$$
In particular, $\{u_0^{\varepsilon}>1\}\supset \{u_1^{\varepsilon}>1\}\neq \emptyset.$
\end{lemma}
{\bf Proof of Lemma 5.2}\\
For $\varepsilon<\varepsilon_0(\lambda),$ let
$$\beta_{\varepsilon}(x,s)=\frac{1}{\varepsilon}\beta\left(\frac{\Min(s,u_0^{\varepsilon}(x))-1}{\varepsilon}\right),\quad B_{\varepsilon}(x,s)=\int_{0}^{s}\beta_{\varepsilon}(x,t)dt,$$
$$\widetilde{g}(x,s)=g_{\varepsilon}(x,(\Min(s,u_0^{\varepsilon}(x))-1)_+),\quad \widetilde{G}_{\varepsilon}(x,s)=\int_{0}^{s}\widetilde{g}(x,t)dt,$$
and set
$$\widetilde{E}_{\varepsilon}(u)=\int_{\Omega}\left[\frac{1}{2}|\nabla_{\mathbb{G}}u|^2+B_{\varepsilon}(x,u)-\lambda \widetilde{G}_{\varepsilon}(x,u)-\beta\varphi_{\beta}(u)\right]dx,\quad u\in W_0^{1,2}(\Omega).$$
Remark that $\widetilde{E}_{\varepsilon}$ is of class $C^1$ and its critical points are the weak solutions of the following problem
\begin{equation}\label{PPP}
  \left\{\begin{array}{ll}  -\mathcal{L} v_0=-\beta_{\varepsilon}(x,u)+\lambda \widetilde{g}(x,u)+\beta \varphi_{\beta}(u) & \quad \mbox{
in }\ \Omega,
\\[0.3cm]u =0 &\quad \mbox{ on }\
  \partial \Omega.
\end{array}
\right.
 \end{equation}
 The elliptic regularity gives that the weak solution of $(\ref{PPP})$ is also a classical solution and the maximum principle gives that $u\leq u_0^{\varepsilon}.$ Hence, $u$ is solution of problem $(\ref{appsol})$ and thus a critical point of $E_{\varepsilon}(u)=\widetilde{E}_{\varepsilon}(u).$\\
 Now, we apply theorem (\ref{Am-Rabin}) to  the functional $\widetilde{E}_{\varepsilon}.$\\
 First, we have $\widetilde{g}_{\varepsilon}(x,s)=\widetilde{g}_{\varepsilon}(x,0)=0$ for $s\leq 1$ and $$\widetilde{g}_{\varepsilon}(x,s)\leq a_0+a_1(\Min(s,u_0^{\varepsilon}(x))-1)_+^{p-1}\leq a_0+a_1(s-1)^{p-1}.$$
 For $s>1,$ and by $g_1),$ we have
 $$\widetilde{G}_{\varepsilon}(x,s)\leq a_0(s-1)_+ +\frac{a_1}{p}(s-1)_+^p $$
 $$\leq (a_0+\frac{a_1}{p})|s|^q$$
 for all $s,$ where $q>2$ if $N=2,$ and $2<q<\frac{2N}{N-2}$ if $N\geq 3.$\\
 So,
 \begin{equation}\label{2etoile}
 \widetilde{E}_{\varepsilon}(u)\geq \int_{\Omega}\left[\frac{1}{2}|\nabla_{\mathbb{G}}u|^2-\lambda(a_0+\frac{a_1}{p})|u|^q-\beta|u|^{1-\delta}\right]dx
 \end{equation}
 $$\geq \frac{1}{2}\|u\|^{2}-\lambda c_1(a_0+\frac{a_1}{p})\|u\|^{q}-\beta c_2\|u\|^{1-\delta}.$$
 By the embedding result $W_0^{1,2}(\Omega)\hookrightarrow L^q(\Omega)$ for $q>2,$ the infimum $m_2(\lambda)$ of the integral $(\ref{2etoile})$ is positive if $\|u\|=r$ when $u\in \partial B_r(0)$ for $r$ sufficiently small when $B_r(0)=\{u\in W_0^{1,2}(\Omega): \|u\|<r \}.$\\
 Taking $r<\|u_0^{\varepsilon}\|$ and applying theorem (\ref{Am-Rabin}), we get a critical point $u_1^{\varepsilon}$ of $\widetilde{E}_{\varepsilon}$ with
 $$\widetilde{E}_{\varepsilon}=\Inf_{\gamma\in \Gamma}\Max_{u\in \gamma([0,1])} \widetilde{E}_{\varepsilon}(u)\geq \Inf_{u\in \partial B_r(0)}\widetilde{E}_{\varepsilon}(u)\geq m_2(\lambda),$$
 where $\Gamma=\{\gamma\in C([0,1],W_0^{1,2}(\Omega)): \gamma(0)=0,\gamma(1)=u_0^{\varepsilon}\}$ is the class of paths joining $0$ and $u_0^{\varepsilon}.$ For the path $\gamma_0(t)=t u_0^{\varepsilon},$ $t\in [0,1],$ we have
 $$\widetilde{E}_{\varepsilon}(\gamma_0(t))\leq \int_{\Omega}\left(\frac{1}{2}|\nabla_{\mathbb{G}}u_0^{\varepsilon}|^2+B_{\varepsilon}(x,u_0^{\varepsilon})\right)dx$$
 since $B_{\varepsilon}(x,s)$ is nondecreasing in $s$ and $\widetilde{G}_{\varepsilon}(x,s)\geq 0$ for all $s$ by $g_2).$\\
 In other part, $$B_{\varepsilon}(x,u_0^{\varepsilon})=\int_{0}^{u_0^{\varepsilon}}\frac{1}{\varepsilon}\beta\left(\frac{t-1}{\varepsilon}\right)dt$$
 $$=B\left(\frac{u_0^{\varepsilon}(x)-1}{\varepsilon}\right)\leq 1.$$
 So, $$\widetilde{E}_{\varepsilon}(u_1^{\varepsilon})\leq \frac{1}{2}\|u_0^{\varepsilon}\|^2+\mathcal{H}(\Omega).$$
 We are now prepared to establish Theorem $1.2$\\\\
 Let $\lambda>\lambda^*$ and choose a sequence $\varepsilon_j\rightarrow 0$ such that $\varepsilon_j<\varepsilon_0(\lambda).$\\
 The lemma $5.1$ gives a minimizer $u_0^{\varepsilon}>0$ of $E_{\varepsilon_j}$ satisfying
 \begin{equation}\label{car}
 E_{\varepsilon_j}(u_0^{\varepsilon_j})\leq m_1(\lambda)+2\lambda \varepsilon_j a_0 \mathcal{H}(\Omega)<0.
 \end{equation}
The lemma $5.2$ gives a second critical point $0<u_1^{\varepsilon_j}<u_0^{\varepsilon_j}$ satisfying
\begin{equation}\label{carcar}
m_2(\lambda)\leq E_{\varepsilon_j}(u_1^{\varepsilon_j})\leq \frac{1}{2}\|u_0^{\varepsilon_j}\|^{2}+\mathcal{H}(\Omega).
\end{equation}
The following step is to verify that $(u_0^{\varepsilon_j})$ and $(u_1^{\varepsilon_j})$ are bounded in $W_0^{1,2}(\Omega)\cap L^{\infty}(\Omega)$ to apply lemma $4.1.$\\
Remark that by $g_1),$
$$G_{\varepsilon}(x,(s-1)_+)\leq a_0(s-1)_++\frac{a_1}{p}(s-1)_+^p$$
$$\leq (a_0+\frac{a_1}{p})|s|^p,\quad\hbox{for all}\quad s.$$
Hence, $$\frac{1}{2}\|u_0^{\varepsilon_j}\|^{2}\leq E_{\varepsilon}(u_0^{\varepsilon})+\lambda (a_0+\frac{a_1}{p})\int_{\Omega}(u_0^{\varepsilon})^pdx.$$
Since $E_{\varepsilon_j}(u_0^{\varepsilon_j})<0$ by $(\ref{car})$ and $p<2,$ then $u_0^{\varepsilon_j}$ is bounded in $W_0^{1,2}(\Omega).$\\
With the same argument and using $(\ref{carcar}),$ we prove that $ u_1^{\varepsilon_j}$ is bounded in $W_0^{1,2}(\Omega).$\\
Also, $$g_{\varepsilon}(x,(s-1)_+)=g_{\varepsilon}(x,0)\quad\hbox{for}\quad s\leq 1 $$
and $$g_{\varepsilon}(x,(s-1)_+)\leq a_0+a_1(s-1)^{p-1}\leq (a_0+a_1)s^{p-1}\quad\hbox{for}\quad s> 1.$$
So, we obtain
$$-\mathcal{L}u_0^{\varepsilon_j}=-\frac{1}{\varepsilon_j}\beta(\frac{u_0^{\varepsilon_j}-1}{\varepsilon_j})+\lambda g_{\varepsilon_j}(x,(u_0^{\varepsilon_j}-1)_+)+\beta(u_0^{\varepsilon_j})^{-\delta}$$
$$\leq \lambda (a_0+a_1)(u_0^{\varepsilon_j})^{p-1}+\beta(u_0^{\varepsilon_j})^{-\delta}.$$
From lemma $5.6$ of \cite{Ghosh}, we conclude that $(u_0^{\varepsilon_j})$ is also bounded in $L^{\infty}(\Omega).$\\
Thus, $(u_1^{\varepsilon_j})$ is bounded because $0<u_1^{\varepsilon_j}\leq u_0^{\varepsilon_j}.$\\
For a relabeled subsequence of $(\varepsilon_j),$ $u_0^{\varepsilon_j}$ and $u_1^{\varepsilon_j}$ converge uniformly  to Lipschitz continuous functions $u_0,u_1\in W_0^{1,2}(\Omega)\cap L^{\infty}(\Omega)$ of problem $(\ref{3}).$\\\\
Furthermore,
\begin{equation}\label{tri}
 E(u_0)\leq \Lim \inf E_{\varepsilon_j}(u_0^{\varepsilon_j})\leq \Lim \sup  E_{\varepsilon_j}(u_0^{\varepsilon_j})\leq E(u_0)+\mathcal{H}\{u_0=1\})
\end{equation}
and
\begin{equation}\label{tritri}
E(u_1)\leq \Lim \inf E_{\varepsilon_j}(u_1^{\varepsilon_j})\leq \Lim \sup  E_{\varepsilon_j}(u_0^{\varepsilon_j})\leq E(u_1)+\mathcal{H}(\{u_1=1\})
\end{equation}
Combining $(\ref{tri})$ and $(\ref{car})$ with $(\ref{neg}),$ we have
$$E(u_0)\leq \Lim \sup E_{\varepsilon_j}(u_0^{\varepsilon_j})\leq m_1(\lambda)\leq E(u_0),$$
hence,
\begin{equation}\label{cir}
E(u_0)=m_1(\lambda)<-\mathcal{H}(\Omega).
\end{equation}
Now, combining $(\ref{tritri})$ and $(\ref{carcar}),$ we obtain
$$E(u_1)+\mathcal{H}(\{u_1=1\})\geq \Lim\inf E_{\varepsilon_j}(u_1^{\varepsilon_j})\geq m_2(\lambda)>0.$$
So, it follows that
\begin{equation}\label{circir}
E(u_1)>-\mathcal{H}(\{u_1=1\})\geq -\mathcal{H}(\Omega).
\end{equation}
Finally, from $(\ref{cir})$ and $(\ref{circir}),$ we conclude that $u_0$ and $u_1$ are distinct and nontrivial solutions of problem $(\ref{3}).$ Also, since $u_1^{\varepsilon_j}\leq u_0^{\varepsilon_j}$ for each $j,$ we have $u_1\leq u_0$ and the sets $\{u_0<1\}\subset \{u_1<1\}$ are connected if $\partial \Omega$ is connected. The sets $\{u_0>1\}\supset \{u_1>1\}$ are nonempty.
\section*{Ethical Approval}
This declaration is not applicable.
\section*{Competing interests}
The authors declare that there are no competing interests between them.
\section*{ Authors' contributions }
The authors declare that they read and approved the final manuscript.
\section*{ Funding}
The authors received no direct funding for this work.
\section*{ Availability of data and materials}
This declaration is not applicable.

\end{document}